\documentclass[11pt]{article}
\usepackage{amsthm}
\usepackage{fullpage}

\usepackage{amsmath}
\usepackage{amsfonts}
\usepackage{amssymb}
\usepackage{latexsym}
\usepackage[dvips]{graphicx}
\usepackage{url}

\theoremstyle{plain}
\newtheorem{thm}{Theorem}
\newtheorem{cor}[thm]{Corollary}
\newtheorem{lem}[thm]{Lemma}

\newcounter{ex_counter}
\newcounter{def_counter}

\theoremstyle{definition}
\newtheorem{ex}[ex_counter]{Example}
\newtheorem{defn}[def_counter]{Definition}

\newcommand{\Q}{\mathbb Q}

\newcommand{\mb}[1]{\mathbf{#1}}
\newcommand{\SW}{\textrm{Sub}} % the set of subwords
\newcommand{\OC}{\textrm{Occ}} % the set of occurrences

\newcommand{\inv}[1]{\textrm{inv}(#1)} %inverse mordent

\title{The critical exponent of the Arshon words}
\author{Dalia Krieger\\
        School of Computer Science \\
        University of Waterloo\\
        Waterloo, ON~N2L 3G1, CANADA\\
        {\tt d2kriege@cs.uwaterloo.ca}}
\begin {document}
\maketitle
\begin{abstract}
Generalizing the results of Thue (for $n = 2$) and of Klepinin and Sukhanov (for $n = 3$), we prove
that for all $n\geq 2$, the critical exponent of the Arshon word of order $n$ is given by
$(3n-2)/(2n-2)$, and this exponent is attained at position 1.
\end{abstract}

\section{Introduction}
In 1935, the Russian mathematician Solomon Efimovich Arshon\footnote{Vilenkin, in his 1991 article
``Formulas on cardboard" \cite{Vilenkin:1991}, says that Arshon was arrested by the Soviet regime
and died in prison, most likely in the late 1930's or early 1940's.} \cite{Arshon:1935,Arshon:1937}
gave an algorithm to construct an infinite cube-free word over 2 letters, and an algorithm to
construct an infinite square-free word over $n$ letters for each $n\geq 3$. The binary word he
constructed turns out to be exactly the celebrated \emph{Thue-Morse word}, $\mb t = 01101001\cdots$
\cite{Thue:1912,Berstel:1995}; the square-free words are now known as the \emph{Arshon words}. For
$n\geq 2$, the Arshon word of order $n$ is denoted by $\mb a_n = a_0a_1a_2\cdots$.

A \emph{square} or a \emph{$2$-power} is a word of the form $xx$, where $x$ is a nonempty word.
Similarly, an \emph{$n$-power} is a word of the form $x^n = xx\cdots x$ ($n$ times). The notion of
integral powers can be generalized to fractional powers. A non-empty finite word $z$ over a finite
alphabet $\Sigma$ is a \emph{fractional power} if it has the form $z = x^ny$, where $x$ is a
non-empty word, $n$ is a positive integer, and $y$ is a prefix of $x$, possibly empty. If $|z| = p$
and $|x| = q$, we say that $z$ is a $(p/q)$-power, or $z = x^{p/q}$. Let $\alpha > 1$ be a real
number. A right-infinite word $\mb w$ over $\Sigma$ is said to be \emph{$\alpha$-power-free} (resp.
$\alpha^+$-power-free), or to \emph{avoid $\alpha$-powers} (resp. $\alpha^+$-powers), if no subword
of it is an $r$-power for any rational $r\geq\alpha$ (resp. $r > \alpha$). Otherwise, \emph{$\mb w$
contains an $\alpha$-power}. The \emph{critical exponent} of $\mb w$, denoted by $E(\mb w)$, is the
supremum of the set of exponents $r\in\Q_{\geq1}$, such that $\mb w$ contains an $r$-power; it may
or may not be attained.

Arshon constructed the words $\mb a_n$, $n\geq 3$, as square-free infinite words. But actually,
these words avoid smaller powers. In 2001, Klepinin and Sukhanov \cite{Klepinin&Sukhanov:2001}
proved that $E(\mb a_3) = 7/4$, and the bound is attained; that is, $\mb a_3$ avoids
$(7/4)^+$-powers. In this paper we generalize the result of Klepinin and Sukhanov, and prove the
following theorem:

\begin{thm}\label{thm:ArshonCE}
Let $n\geq 2$, and let $\mb a_n = a_0a_1a_2\cdots$ be the Arshon word of order $n$. Then the
critical exponent of $\mb a_n$ is given by $E(\mb a_n) = (3n-2)/(2n-2)$, and $E(\mb a_n)$ is
attained by a subword beginning at position 1.
\end{thm}

\section{Definitions and notation}
\subsection{Definition of the Arshon words} Let $\Sigma_n = \{0,1,\ldots,n-1\}$ be an alphabet of
size $n$, $n\geq 3$. In what follows, we use the notation $a\pm1$, where $a\in\Sigma_n$, to denote
the next or previous letter in lexicographic order, and similarly we use the notation $a+b$, $a-b$,
where $a,b\in\Sigma_n$; all sums of letters are taken modulo $n$.

Define two morphisms over $\Sigma_n$ as follows:
$$
\begin{array}{llll}
\varphi_{e,n}(a) &=&  a(a+1)\cdots (n-1)01\cdots (a-2)(a-1), &  a = 0,1,\cdots,n-1;\vspace{2mm}\\
\varphi_{o,n}(a) &=&  (a-1)(a-2)\cdots 10(n-1)\cdots (a+1)a, & a = 0,1,\cdots,n-1.
\end{array}
$$
The letters `$e$' and `$o$' stand for ``even" and ``odd", respectively.  Both $\varphi_{e,n}$ and
$\varphi_{o,n}$ are \emph{$n$-uniform} (that is, $|\varphi_{e,n}(a)| = n$ for all $a\in\Sigma_n$,
and similarly for $\varphi_{o,n}$) and \emph{marked} (that is, $\varphi_{e,n}(a)$ and
$\varphi_{e,n}(b)$ have no common prefix or suffix for all $a\neq b\in\Sigma_n$, and similarly for
$\varphi_{o,n}$). The Arshon word of order $n$ can be generated by alternately iterating
$\varphi_{e,n}$ and $\varphi_{o,n}$: define an operator $\varphi_n:\Sigma_n^*\rightarrow\Sigma_n^*$
by
\begin{equation}\label{eq:ArshonOp}
\varphi_n(a_i) =
  \left\{
  \begin{array}{ll}
    \varphi_{e,n}(a_i), & \textrm{ if } i\textrm{ is even};\\
    \varphi_{o,n}(a_i), & \textrm{ if } i\textrm{ is odd}.\\
  \end{array}\right.
\end{equation}
That is, if $u = a_0a_1\cdots a_m\in\Sigma_n^*$, then $\varphi_n(u) =
\varphi_{e,n}(a_0)\varphi_{o,n}(a_1)\varphi_{e,n}(a_2)\varphi_{o,n}(a_3)\cdots$. The Arshon word of
order $n$ is given by
\begin{equation}\label{eq:ArshonWord}
\mb a_n = \lim_{k\rightarrow\infty}\varphi_n^k(0).
\end{equation}
Note that $\varphi_n^k(0)$ is a prefix of $\varphi_n^{k+1}(0)$ for all $k \geq 0$, and the limit is
well defined.

\begin{ex}\label{ex:a3}
For $n=3$, the even and odd Arshon morphisms are given by
$$
 \varphi_{e,3}:\left\{
  \begin{array}{lll}
    0 & \rightarrow & 012 \\
    1 & \rightarrow & 120 \\
    2 & \rightarrow & 201\\
  \end{array} \right.\;,\;\;\;
\varphi_{o,3}:\left\{
  \begin{array}{lll}
    0 & \rightarrow & 210 \\
    1 & \rightarrow & 021 \\
    2 & \rightarrow & 102\\
  \end{array}\right.\;,\;\;\;
$$
and the Arshon word of order 3 is given by
$$
\mb a_3 = \lim_{k\rightarrow\infty}\varphi_3^k(0) =
    \underbrace{012}_{\varphi_{e,3}(0)}
    \underbrace{021}_{\varphi_{o,3}(1)}
    \underbrace{201}_{\varphi_{e,3}(2)}
    \underbrace{210}_{\varphi_{o,3}(0)}\cdots.
$$
\end{ex}

It is not difficult to see that when $n$ is even, the $i$'th letter of $\mb a_n$ is even if and
only if $i$ is an even position (for a formal proof, see S\'e\'ebold
\cite{Seebold:2002,Seebold:2003}). Therefore, when $n$ is even, the map $\varphi_n$ becomes a
morphism, denoted by $\alpha_n$:
\begin{equation}\label{eq:ArshonMorph}
\alpha_n(a) =
  \left\{
  \begin{array}{ll}
    \varphi_{e,n}(a), & \textrm{ if } a\textrm{ is even};\\
    \varphi_{o,n}(a), & \textrm{ if } a\textrm{ is odd}.\\
  \end{array}\right.
\end{equation}
When $n$ is odd no such partition exists, and indeed, $\mb a_n$ cannot be generated by iterating a
morphism. This fact was proved for $\mb a_3$ by Berstel \cite{Berstel:1979} and Kitaev
\cite{Kitaev:2000,Kitaev:2003}, and for any odd $n$ by Currie \cite{Currie:2002}.

\subsection{Subwords and occurrences}
An \emph{occurrence} of a subword within $\mb a_n$ is a triple $(z,i,j)$, where $z$ is a subword of
$\mb a_n$, $0\leq i\leq j$, and $a_i\cdots a_j = z$. In other words, $z$ occurs in $\mb a_n$ at
positions $i,\cdots,j$. We usually refer to an occurrence $(z,i,j)$ as $z = a_i\cdots a_j$. The set
of all subwords of $\mb a_n$ is denoted by $\SW(\mb a_n)$. The set of all occurrences of subwords
within $\mb a_n$ is denoted by $\OC(\mb a_n)$. An occurrence $(z,i,j)$ \emph{contains} an
occurrence $(z',i',j')$ if $i\leq i'$ and $j\geq j'$.

A subword $v$ of $\mb a_n$ admits an \emph{interpretation} by $\varphi_n$ if there exists a subword
$v' = v_0v_1\cdots v_kv_{k+1}$ of $\mb a_n$, $v_i\in \Sigma_n$, such that $v =
y_0\varphi_n(v_1\cdots v_k)x_{k+1}$, where $y_0$ is a suffix of $\varphi_n(v_0)$ and $x_{k+1}$ is a
prefix of $\varphi_n(v_{k+1})$. The word $v'$ is called an \emph{ancestor} of $v$.

For an occurrence $z\in\OC(\mb a_n)$, we denote by $\inv{z}$ the inverse image of $z$ under
$\varphi_n$. That is, $\inv{z}$ is the shortest occurrence $z'\in\OC(\mb a_n)$ such that
$\varphi_n(z')$ contains $z$. Note that the word (rather than occurrence) $\inv{z}$ is an ancestor
of the word $z$, but not necessarily a unique one.

Following Currie \cite{Currie:2002}, we refer to the decomposition of $\mb a_n$ into images under
$\varphi_n$ as the \emph{$\varphi$-decomposition}, and to the images of the letters as
\emph{$\varphi$-blocks}. We denote the borderline between two consecutive $\varphi$-blocks by
`$|$'; e.g., $i|j$ means that $i$ is the last letter of a block and $j$ is the first letter of the
following block. If $z = a_i\cdots a_j\in\OC(\mb a_n)$ begins at an even position we write $z =
a_i^{(e)}a_{i+1}^{(o)}a_{i+2}^{(e)}\cdots$, and similarly for an occurrence that begins at an odd
position.

\section{General properties of the Arshon words}\label{sect:general}

\begin{lem}\label{lem:ArshonLowerBound}
For all $n\geq 2$, $\mb a_n$ contains a $(3n-2)/(2n-2)$-power beginning at position 1.
\end{lem}
\begin{proof}
For $n = 2$, $\mb a_2 = \mb t = 0110\cdots$, which contains the 2-power $11$ at position 1. For
$n\geq 3$, $\mb a_n$ begins with
\begin{multline}\nonumber
\varphi_{e,n}(0)\varphi_{o,n}(1)\varphi_{e,n}(2) = 012\cdots(n-1)|0(n-1)\cdots21|2\cdots(n-1)01 =\\
0\; (12\cdots(n-1)0(n-1)\cdots2)^{(3n-2)/(2n-2)} \; 1.
\end{multline}
\end{proof}

\begin{ex}\label{ex:ArshonLowerBound}
$$
\begin{array}{lllllllll}
\mb a_3 &\;& = &\;& 012|021|201|\cdots       &\;& = &\;& 0\;(1202)^{7/4}\;1\cdots, \vspace{1mm}\\
\mb a_4 &\;& = &\;& 0123|0321|2301|\cdots    &\;& = &\;& 0\;(123032)^{10/6}\;1\cdots, \vspace{1mm}\\
\mb a_5 &\;& = &\;& 01234|04321|23401|\cdots &\;& = &\;& 0\;(12340432)^{13/8}\;1\cdots. \vspace{1mm}\\
\end{array}
$$
\end{ex}

\begin{cor}\label{cor:ArshonLowerBound}
The critical exponent of $\mb a_n$ satisfies $(3n-2)/(2n-2) \leq E(\mb a_n) \leq 2$ for all $n\geq
2$.
\end{cor}
\begin{proof}
For $n = 2$, it is well known that $E(\mb a_n) = E(\mb t) = 2$ \cite{Thue:1912,Berstel:1995}. For
$n\geq 3$, we know by Arshon \cite{Arshon:1935,Arshon:1937} that $\mb a_n$ is square-free, and so
$E(\mb a_n) \leq 2$. The lower bound follows from Lemma~\ref{lem:ArshonLowerBound}. \end{proof}

\begin{lem}\label{lem:ArshonStruct}
Let $n\geq 3$, and let $i,j\in\Sigma_n$.
\begin{enumerate}
  \item If $ij\in\OC(\mb a_n)$, then $j = i\pm1$.
  \item The borderline between two consecutive $\varphi$-blocks has the form $i|ji$ or $ij|i$.
  Moreover, a word of the form $iji$ can occur only at a borderline.
\end{enumerate}
\end{lem}
\begin{proof}
If $ij$ occurs within a $\varphi$-block, then $j = i\pm1$ by definition of $\varphi_n$. Suppose $i$
is the last letter of a $\varphi$-block and $j$ is the first letter of the next $\varphi$-block,
and let $kl = \inv{ij}$. Assume $j\neq i\pm1$, and suppose further that $ij$ is the first pair that
satisfies this inequality. Then $l= k\pm1$, and so there are four cases:
$$
\begin{array}{lllllllll}
\varphi_n(kl) &=& \varphi_{e,n}(k)\varphi_{o,n}(k+1), & \varphi_n(kl) &=& \varphi_{o,n}(k)\varphi_{e,n}(k+1),\vspace{2mm}\\
\varphi_n(kl) &=& \varphi_{e,n}(k)\varphi_{o,n}(k-1), & \varphi_n(kl) &=&
\varphi_{o,n}(k)\varphi_{e,n}(k-1).
\end{array}
$$
But it is easy to check that for all the cases above, $j = i\pm1$, a contradiction.

For the second assertion, observe that by definition of $\varphi_n$, a $\varphi$-block is either
strictly increasing or strictly decreasing, and two consecutive blocks have alternating directions.
By the above, a change of direction can have only the form $i|ji$ or $ij|i$.
\end{proof}

\begin{defn}[Currie \cite{Currie:2002}]\label{def:mordent}
A \emph{mordent} is a word of the form $iji$, where $i,j\in\Sigma_n$ and $j = i\pm1$. Two
consecutive mordents occurring in $\mb a_n$ are either \emph{near mordents}, \emph{far mordents},
or \emph{neutral mordents}, according to the position of the borderlines:
\begin{eqnarray}
\nonumber i|ji \;u\; kl|k &=& \textrm{near mordents, }|u| = n-4;\\
\nonumber ij|i \;u\; k|lk &=& \textrm{far mordents, }|u| = n-2;\\
\nonumber i|ji \;u\; k|lk &=& \textrm{neutral mordents, }|u| = n-3;\\
\nonumber ij|i \;u\; kl|k &=& \textrm{neutral mordents, }|u| = n-3.
\end{eqnarray}
Note that for $n = 3$, near mordents are overlapping: $ \mb a_3 =
01\emph{2}|\emph{0}\emph{\textbf{2}}\emph{1}|\emph{2}01|\cdots$.
\end{defn}

Since $\mb a_n$ is square-free, a $p/q$-power occurring in $\mb a_n$ has the form $xyx$, where $q =
|xy|$, $p = |xyx|$, and both $x$, $y$ are nonempty.

\begin{defn} \label{def:stretch}
Let $z = a_i\cdots a_j \in \OC(\mb a_n)$ be a $p/q$-power. We say that $z$ is
\emph{left-stretchable} (resp. \emph{right-stretchable}) if the $q$-period of $z$ can be stretched
left (resp. right), i.e., if $a_{i-1} =  a_{i+q-1}$ (resp. $a_{j+1} = a_{j-q+1}$). If  the
$q$-period of $z$ can be stretched neither left nor right, we say that $z$ is an
\emph{unstretchable} $p/q$-power.
\end{defn}

Since the critical exponent is a supremum, it is enough to consider unstretchable powers when
computing it.

\begin{lem}\label{lem:ArshonCE}
Let $n\geq 3$. Let $z = xyx = (xy)^{p/q}\in\OC(\mb a_n)$ be an unstretchable power such that
$|x|\leq n$ and $x$ contains no mordents. Then $p/q \leq (3n~-~2)/(2n~-~2)$.
\end{lem}
\begin{proof}
Since $|x| \leq n$, it is enough to consider $y$ such that $|y|\leq n-2$, for otherwise we would
get that $p/q < (3n-2)/(2n-2)$. Therefore, $|xy| = q \leq 2n-2$ and $|z| \leq 3n-2$. We get that
$xy$ is contained in at most 3 consecutive $\varphi$-blocks and $z$ is contained in at most 4
consecutive $\varphi$-blocks. Suppose $z$ is not contained in 3 consecutive $\varphi$-blocks. Let
$B_0B_1B_2B_3$ be the blocks containing $z$, and assume that $B_0$ is even (the other case is
similar). Since $|x|\leq n$, necessarily $xy$ begins in $B_0$ and ends in $B_2$. Since $x$ contains
no mordents, $x$ has to start at the last letter of $B_0$: otherwise, we would get that $x$ cannot
extend beyond the first letter of $B_1$, and since $|y|\leq n-2$, we would get that $z$ is
contained in 3 $\varphi$-blocks. Therefore, the letters of $x$ are decreasing. Now, since $|xy|\leq
2n-2$, the second occurrence of $x$ begins at least 3 letters from the end of $B_2$. Since $B_2$ is
an even block, we get a contradiction if $|x| > 1$. But if $|x| = 1$ then $z$ is contained in
$B_0B_1B_2$. We can assume therefore that $z$ is contained in 3 consecutive $\varphi$-blocks,
$B_0B_1B_2$. We assume that $B_0$ is even (the other case is symmetric).

If $xy$ is contained in one block then, because $B_0,B_2$ are even and $B_1$ is odd, necessarily
$|x| = 1$, and so $p/q \leq 3/2 < (3n-2)/(2n-2)$. If $xy$ begins in $B_0$ and ends in $B_2$, then,
since $|y|\leq n-2$, the first $x$ occurrence has to end at the third letter of $B_1$ or later.
Since $x$ contains no mordents, this implies that $xy$ begins at the last letter of $B_0$ and the
letters of $x$ are decreasing. Since $B_2$ is even, again $|x| = 1$.

Assume $xy$ begins in $B_0$ and ends in $B_1$. Again, because $B_0$ is even and $B_1$ is odd, in
order for $x$ to contain more than one letter the second $x$ occurrence has to start either at the
last letter of $B_1$, or at the first letter of $B_2$.

Let $B_0 = \varphi_{e,n}(i)$. Then there are four cases for $B_1, B_2$:
\begin{enumerate}
  \item $B_1 = \varphi_{o,n}(i+1)$, $B_2 = \varphi_{e,n}(i)$;
  \item $B_1 = \varphi_{o,n}(i-1)$, $B_2 = \varphi_{e,n}(i)$;
  \item $B_1 = \varphi_{o,n}(i+1)$, $B_2 = \varphi_{e,n}(i+2)$;
  \item $B_1 = \varphi_{o,n}(i-1)$, $B_2 = \varphi_{e,n}(i-2)$.
\end{enumerate}
We now check what the maximal possible exponent is in each of these cases. Without loss of
generality, we can assume $i = 0$. We use the notation $z = xyx'$, where $x'$ is the second
occurrence of $x$ in $z$.

\begin{description}
  \item[Case 1:] $B_0B_1B_2 = |01\cdots(n-1)|0(n-1)\cdots 1|01\cdots(n-1)|\;$.

If $x'$ starts at the last letter of $B_1$ then $|x| = 1$, since $10$ does not occur anywhere
before. If $x'$ starts at the first letter of $B_2$, the only possible power is the $3n/2n$-power
$B_0B_1B_0$, which contradicts the hypothesis $|y|\leq n-2$.

  \item[Case 2:] $B_0B_1B_2  = |01\cdots(n-1)|(n-2)(n-3)\cdots0(n-1)|01\cdots(n-1)|\;$.

By the same argument, either $|x| = 1$ or $z$ is a $3n/2n$-power.

  \item[Case 3:] $B_0B_1B_2 = |01\cdots(n-1)|0(n-1)\cdots 1|23\cdots(n-1)01|\;$.

If $x'$ starts at the last letter of $B_1$, we get the $(3n-2)/(2n-2)$-power described in
Lemma~\ref{lem:ArshonLowerBound}. If $x'$ starts at the first letter of $B_2$, then $x$ has to
start at the $2$ in $B_0$. But then the power is left-stretchable, to the $(3n-2)/(2n-2)$-power
described above.

  \item[Case 4:] $B_0B_1B_2 = |01\cdots(n-1)|(n-2)(n-3)\cdots0(n-1)|(n-2)(n-1)0\cdots(n-3)|\;$.

If $x'$ starts at the last letter of $B_1$, then $x$ has to start at the last letter of $B_0$. But
then $|x| =2$, since $(n-1)\neq (n-3)$. We get that $z$ is an $(n+2)/n$-power, and $(n+2)/n <
(3n-2)/(2n-2)$ for all $n\geq 3$. If $x'$ starts at the first letter of $B_2$, then $x$ has to
start at the second last letter of $B_0$. Again, $|x| = 2$, and $z$ is an $(n+4)/(n+2)$-power,
where $(n+4)/(n+2) < (3n-2)/(2n-2)$ for all $n \geq 2$.
\end{description}
\end{proof}

In what follows, we will show that in order to compute $E(\mb a_n)$, it is enough to consider
powers $xyx$ such that $|x| \leq n$ and $x$ contains no mordents.

\begin{defn}\label{def:sybchpoint}
Let $z$ be a subword of $\mb a_n$. We say that $(z_1,z_2)$ is a \emph{synchronization point} of $z$
under $\varphi_n$ if $z = z_1z_2$, and whenever $\varphi_n(u) = v_1zv_2$ for some
$u,v_1,v_2\in\SW(\mb a_n)$, we have $u = u_1u_2$, $\varphi_n(u_1) =  v_1z_1$, and $\varphi_n(u_2) =
z_2v_2$. That is, $z_1|z_2$ is always a borderline in the $\varphi$-decomposition of $z$,
regardless of the position in $\mb a_n$ where $z$ occurs. We say that a subword $z\in\SW(\mb a_n)$
is \emph{synchronized} if it can be decomposed unambiguously under $\varphi_n$, in which case it
has a unique ancestor.
\end{defn}

\begin{lem}\label{lem:ArshonSynch}
If $z\in\SW(\mb a_n)$ has a synchronization point then $z$ is synchronized.
\end{lem}
\begin{proof}
Suppose $z$ has a synchronization point, $z = u|v$. If $|u| = |v| = 1$ then $z$ cannot have a
synchronization point at $u|v$, since $uv$ occurs either in $\varphi_{e,n}(u)$ or in
$\varphi_{o,n}(u+1)$. Therefore, at least one of $u$, $v$ has length $> 1$. Suppose $|u| > 1$. If
the last two characters of $u$ are increasing, we know that an even $\varphi$-block ends at $u$ and
an odd $\varphi$-block starts at $v$, and vice versa if the last two characters of $u$ are
decreasing. Since both $\varphi_{e,n}$ and $\varphi_{o,n}$ are uniform marked morphisms, and since
we know $\varphi$-blocks alternate between even and odd, we can infer $\inv{z}$ unambiguously from
$u|v$. \end{proof}

\begin{lem}\label{lem:ArshonInvImage}
Let $n\geq3$, and let $z = xyx = (xy)^{p/q}\in\OC(\mb a_n)$ be an unstretchable $p/q$-power, such
that $x$ has a synchronization point. Then there exists an $r/s$-power $z'\in\OC(\mb a_n)$, such
that $p = nr$, $q = ns$, and $z = \varphi_n(z')$.
\end{lem}
\begin{proof}
Since $x$ has a synchronization point, it has a unique decomposition under $\varphi_n$. Suppose $x$
does not begin at a borderline of $\varphi$-blocks. Then $x = t|w$, where $t$ is a nonempty suffix
of a $\varphi$-block, and $z = t|wyt|w$. But since the interpretation is unique, both occurrences
of $t$ must be preceded by a word $s$, such that $st$ is a $\varphi$-block. Thus $z$ can be
stretched by $s$ to the left, a contradiction. Therefore, $x$ begins at a borderline, and so $y$
ends at a borderline. For the same reason, $x$ must end at a borderline, and so $y$ must begin at a
borderline. We get that both $x$ and $y$ have an exact decomposition into $\varphi$-blocks, and
this decomposition is unique. In particular, both occurrences of $x$ have the same inverse image
under $\varphi_n$. Let $k,l$ be the number of $\varphi$-blocks composing $x,y$, respectively. Then
$p = n(2k+l)$, $q = n(k+l)$, and $\varphi_n^{-1}(z) =
\varphi_n^{-1}(x)\varphi_n^{-1}(y)\varphi_n^{-1}(x)$ is a $(2k+l)/(k+l)$-power. \end{proof}

\begin{cor}\label{cor:noSynch}
To compute $E(\mb a_n)$, it is enough to consider powers $z = xyx$ such that $x$ has no
synchronization points.
\end{cor}

\section{Arshon words of even order}\label{sect:ArshonEven}
To illustrate the power structure in Arshon words of even order, consider $\mb a_4$:
$$
\begin{array}{ccccccccc}
  &&&&            \varphi_{e,4} &\;\;& \varphi_{o,4} &\;\;& \alpha_4   \vspace{2mm}\\
0 &\;\;& \rightarrow   &\;\;& 0123 &\;\;& 3210 &\;\;& 0123\\
1 &\;\;& \rightarrow   &\;\;& 1230 &\;\;& 0321 &\;\;& 0321\\
2 &\;\;& \rightarrow   &\;\;& 2301 &\;\;& 1032 &\;\;& 2301\\
3 &\;\;& \rightarrow   &\;\;& 3012 &\;\;& 2103 &\;\;& 2103\\
\end{array}
$$
$$
\mb a_4 = 0123|0321|2301|2103|0123|2103|2301|0321|2301|2103|0123|0321|2301|0321|\cdots
$$

\begin{lem}\label{lem:ArshonEvenSynch}
Let $n\geq4$ be even, and let $x\in\SW(\mb a_n)$ be a subword that has no synchronization point.
Then $|x| \leq n$ and $x$ contains no mordents.
\end{lem}
\begin{proof}
In general, a mordent $iji$ can admit two possible borderlines: $ij|i$ or $i|ji$. However, if $n$
is even, all images under $\alpha_n$ begin with an even letter and end with an odd letter; images
of odd letters under $\varphi_{e,n}$ and images of even letters under $\varphi_{o,n}$ are never
manifested. Therefore, every mordent admits exactly one interpretation: if $i$ is even and $j$ is
odd the interpretation has to be $ij|i$, and vice versa for odd $i$. Thus, if $x$ has no
synchronization point it contains no mordents.

Suppose $x$ contains no mordents. Then $|x| \leq n+2$, and the letters of $x$ are either increasing
or decreasing. Assume they are increasing. If $|x| = n+2$ then $x$ has exactly one interpretation,
$x = i|(i+1)\cdots (i-1)i|(i+1)$, or else we would get that $\mb a_n$ contains two consecutive even
blocks. If $|x| = n+1$ then a priori $x$ has two possible interpretations: $x = i|(i+1)\cdots
(i-1)i|$ or $x = |i(i+1)\cdots (i-1)|i$. However, the first case is possible if and only if $i$ is
odd, since for an even $n$ no $\varphi$-block ends with an even letter. Similarly, the second case
is possible if and only if $i$ is even. \end{proof}

Lemma~\ref{lem:ArshonEvenSynch}, together with Corollary~\ref{cor:noSynch} and
Lemma~\ref{lem:ArshonCE}, completes the proof of Theorem~\ref{thm:ArshonCE} for all even $n\geq 4$.

\section{Arshon words of odd order}\label{sect:ArshonOdd}

To illustrate the power structure in Arshon words of odd order, consider $\mb a_5$:
$$
\begin{array}{ccccccccc}
  &&&&                \varphi_{e,5} &\;\;& \varphi_{o,5} \vspace{2mm}\\
0 &\;\;& \rightarrow   &\;\;& 01234 &\;\;& 43210 \\
1 &\;\;& \rightarrow   &\;\;& 12340 &\;\;& 04321 \\
2 &\;\;& \rightarrow   &\;\;& 23401 &\;\;& 10432 \\
3 &\;\;& \rightarrow   &\;\;& 34012 &\;\;& 21043 \\
4 &\;\;& \rightarrow   &\;\;& 40123 &\;\;& 32104 \\
\end{array}
$$
$$
\mb a_5 = 01234|04321|23401|21043|40123|43210|40123|21043|23401|04321|23401|21043|\cdots
$$

\begin{lem}\label{lem:ArshonOddSynch}
Let $n\geq3$ be odd. Then every subword $z\in\SW(\mb a_n)$ with $|z|\geq 3n$ has a unique
interpretation under $\varphi_n$.
\end{lem}
\begin{proof}
Consider a subword that contains a pair of consecutive mordents, $z = iji\;u\;klk$. If $|u| = n-4$
(that is, these are near mordents), then $z$ contains two synchronization points, $z =
i|ji\;u\;kl|k$: otherwise, we get a $\varphi$-block that contains a repeated letter, a
contradiction. Similarly, if $|u| = n-2$ (a pair of far mordents), $z$ contains the synchronization
points $z = ij|i\;u\;k|lk$. To illustrate, consider $\mb a_5$: let $z = a_4\cdots a_{10} =
404\;3\;212$. A borderline $40|4$ implies that $43212$ is a $\varphi$-block, a contradiction; a
borderline $2|12$ implies that $40432$ is a $\varphi$-block, again a contradiction. Now let $z =
a_7\cdots a_{16} = 212\;340\;121$. A borderline $2|12$ implies that $121$ is a prefix of a
$\varphi$-block, while a borderline $12|1$ implies that $212$ is a suffix of a $\varphi$-block.
Again, we get a contradiction.

If $|u| = n-3$ (neutral mordents), then $z$ has two possible interpretations, either $z =
i|ji\;u\;k|lk$ or $z = ij|i\;u\;kl|k$. However, by Currie \cite{Currie:2002}, $\mb a_n$ does not
contain two consecutive pairs of neutral mordents: out of three consecutive mordents, at least one
of the pairs is either near or far. (It is also easy to see that this is the case by a simple
inverse image analysis: an occurrence of the form $ij|i\;u\;kl|k\;v\;rs|r$ or
$i|ji\;u\;k|lk\;v\;r|sr$ implies that $\mb a_n$ contains a square of the form $abab$,
$a,b\in\Sigma_n$, a contradiction: by Arshon, $\mb a_n$ is square-free.)

Let $z\in\OC(\mb a_n)$ satisfy $|z| = 3n$. If $z$ contains a pair of near or far mordents, then $z$
has a unique ancestor. Otherwise, $z$ contains a pair of neutral mordents, $iji\;u\;klk$, where
$|u| = n-3$, and there are two possible interpretations: $i|ji\;u\;k|lk$ or $ij|i\;u\;kl|k$. Let
$i'j'i'$ be the mordent on the left of $iji$, and let $k'l'k'$ be the mordent on the right of
$klk$. Since no two consecutive neutral mordents occur, $i'j'i'$ and $k'l'k'$ must form near or far
mordents with $iji$ and $klk$.

If the interpretation is $i|ji\;u\;k|lk$, then $k'l'k'$ forms a near pair with $klk$, while
$i'j'i'$ forms a far pair with $iji$. Therefore, $k'l'k'$ is $n-4$ letters away from $klk$, while
$i'j'i'$ is $n-2$ letters away from $iji$. By assumption, $z$ does not contain a near pair or a far
pair, therefore $z$ can contain at most $n-2$ letters to the right of $klk$, and at most $n$
letters to the left of $iji$. Since $|z| = 3n$, this means that either $z =
j'|i'\;x\;i|ji\;u\;k|lk\;v\;k'$ or $z = |i'\;x\;i|ji\;u\;k|lk\;v\;k'l'|$, where $|x| = n-2$ and
$|v| = n-4$. Similarly if the interpretation is $ij|i\;u\;kl|k$, then either $z =
|j'i'\;v\;ij|i\;u\;kl|k\;x\;k'|$ or $z = i'\;v\;ij|i\;u\;kl|k\;x\;k'|l'$, where $|x| = n-2$ and
$|v| = n-4$. In any case, $z$ contains enough letters to determine if the far mordent is on the
left or on the right, and the interpretation is unique.

\end{proof}

\begin{ex}\label{ex:ArshonOddSynch}
For $n = 5$, the occurrence $z = a_{21}\cdots a_{34} = 01234321040123$, of length $3n-1 = 14$, has
two possible interpretations under $\varphi_5$, as illustrated in Fig.~\ref{fig:arshon5}.
\begin{figure}[h]
\begin{center}
  \includegraphics [width=0.4\textwidth]{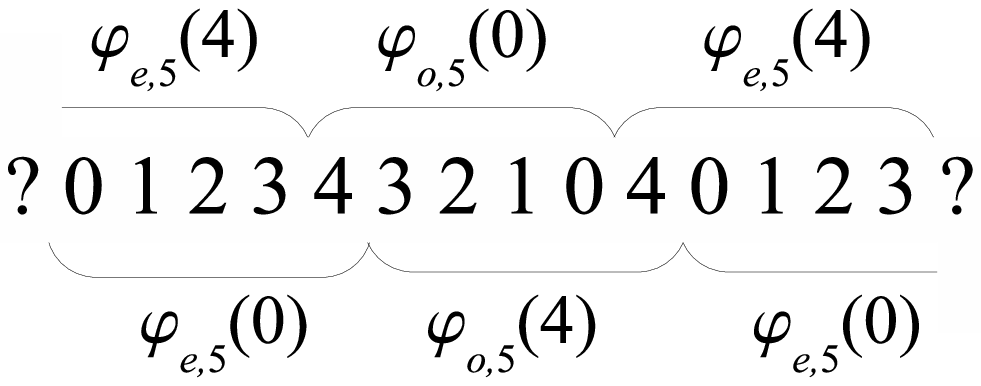}
  \caption{Two interpretations under $\varphi_5$.}
  \label{fig:arshon5}
\end{center}
\end{figure}
However, if either of the left or right question marks is known, the ambiguity is solved: the top
interpretation is valid if and only if the left question mark equals 4 (so as to complete the
$\varphi$-block) and the right question mark equals $2$ (so as to complete the near mordent). The
bottom interpretation is valid if and only if the left question mark equals 1 (so as to complete
the near mordent) and the right question mark equals $4$ (so as to complete the $\varphi$-block).
\end{ex}

\textbf{Note:} Lemma~\ref{lem:ArshonOddSynch} is an improvement of a similar lemma of Currie
\cite{Currie:2002}, who proved that every occurrence of length $3n+3$ or more has a unique
interpretation.

\begin{cor}\label{cor:ArshonOddCE}
The critical exponent of an odd Arshon word is the largest exponent of powers of the form $z =
xyx$, such that $|x| < 3n$.
\end{cor}

To compute $E(\mb a_n)$ we need to consider subwords of the form $xyx$, with $x$ unsynchronized.
Moreover, the two occurrences of $x$ should have different interpretations, or else we could take
an inverse image under $\varphi_n$. For a fixed $n$, it would suffice to run a computer check on a
finite number of subwords of $\mb a_n$; this is exactly the technique Klepinin and Sukhanov
employed in \cite{Klepinin&Sukhanov:2001}. For a general $n$, we need a more careful analysis.
%Since the result for $n = 3$ is known, we concentrate on $n\geq 5$. This will simplify some of the
%arguments.

\begin{lem}\label{lem:mordentInv}
Let $n\geq 3$, $n$ odd. For all mordents in $\mb a_n$,
\begin{enumerate}
  \item $\inv{i(i+1)i} = (i+2)^{(e)}(i+1)^{(o)}$ or $\inv{i(i+1)i} = (i+1)^{(e)}(i+2)^{(o)}$;
  \item $\inv{i(i-1)i} = (i-1)^{(o)}(i)^{(e)}$ or $\inv{i(i-1)i} = (i)^{(o)}(i-1)^{(e)}$.
\end{enumerate}
\end{lem}
\begin{proof}
A mordent $iji$ can admit two possible borderlines: $ij|i$ or $i|ji$. Consider the mordent
$i(i+1)i$. If the borderline is $i(i+1)|i$, then $i(i+1)$ is a suffix of an increasing
$\varphi$-block, and so the block must be an image under $\varphi_{e,n}$. By definition of
$\varphi_{e,n}$, $i(i+1)$ is the suffix of $\varphi_{e,n}(i+2)$. Since even and odd blocks
alternate, the next block must be an image under $\varphi_{o,n}$, and by definition, $i$ is the
prefix of $\varphi_{o,n}(i+1)$.

If the borderline is $i|(i+1)i$, then $(i+1)i$ is the prefix of a decreasing $\varphi$-block, and
by similar considerations this block is $\varphi_{o,n}(i+2)$, while the previous block is
$\varphi_{e,n}(i+1)$. The assertion for $i(i-1)i$ is proved similarly. \end{proof}

\begin{lem}\label{lem:dist}
Let $n\geq 3$, $n$ odd, and let $z\in\OC(\mb a_n)$.
\begin{enumerate}
  \item If $z = i^{(e)}u i^{(o)}$ or $z = i^{(o)}u i^{(e)}$ for some $i\in\Sigma_n$, then $|u| \geq
  n-1$;
  \item If $z = i^{(e)}u (i\pm1)^{(e)}$ or $z = i^{(o)}u (i\pm1)^{(o)}$ for some $i\in\Sigma_n$, then $|u| \geq
  n-2$.
\end{enumerate}
\end{lem}
\begin{proof}
Let $z = i^{(e)}u i^{(o)}$, and suppose $i^{(o)}$ does not occur in $u$ (otherwise, if $u =
u'i^{(o)}u''$, take $z = i^{(e)}u' i^{(o)}$). If $|u| < n-1$ then $z$ must contain a mordent in
order for $i$ to be repeated. But then the two occurrences of $i$ have the same parity, a
contradiction. The rest of the cases are proved similarly. \end{proof}

\begin{lem}\label{lem:noMordent}
Let $n \geq 3$ be odd, and let $z = xyx = (xy)^{p/q}\in\OC(\mb a_n)$ be an unstretchable power,
such that $x$ is unsynchronized and contains a mordent. Then $p/q < E(\mb a_n)$.
\end{lem}
\begin{proof}
%Since $x$ contains one mordent, $|x| \leq 2n+2$.
Suppose $x$ contains the mordent $i(i+1)i$ (the case of $i(i-1)i$ is symmetric). Then the two
occurrences of the mordent have different interpretations, else we could take an inverse image
under $\varphi_n$ and get a power with the same exponent. By Lemma~\ref{lem:mordentInv}, there are
two different cases, according to which interpretation comes first:
\begin{eqnarray}
\nonumber \overbrace{\cdots i(i+1)}^{\varphi_{e,n}(i+2)}|\overbrace{i(i-1)\cdots
(i+2)(i+1)}^{\varphi_{o,n}(i+1)}|&\overbrace{\cdots\cdots\cdots}^{n-1
\;\;\varphi-\textrm{blocks}}&|\overbrace{(i+1)(i+2)\cdots
(i-1)i}^{\varphi_{e,n}(i+1)}|\overbrace{(i+1)i\cdots}^{\varphi_{o,n}(i+2)}\\
\nonumber&&\\
\nonumber\overbrace{\cdots i}^{\varphi_{e,n}(i+1)}|\overbrace{(i+1)i\cdots
(i+3)(i+2)}^{\varphi_{o,n}(i+2)}|&\overbrace{\cdots\cdots\cdots}^{n-1
\;\;\varphi-\textrm{blocks}}&|\overbrace{(i+2)(i+3)\cdots
i(i+1)}^{\varphi_{e,n}(i+2)}|\overbrace{i(i-1)\cdots}^{\varphi_{o,n}(i+1)}
\end{eqnarray}
By Lemma~\ref{lem:dist}, in both cases there must be at least $n-1$ additional $\varphi$-blocks
between the blocks containing the two $i(i+1)i$ occurrences. Therefore, in both cases $q\geq
n^2+n-1$ (note that $q$ is the length of the period, and can be measured from the beginning of
$i(i+1)i$ in the first $x$ to just before $i(i+1)i$ in the second $x$). Now, $x$ is unsynchronized,
and so by Lemma~\ref{lem:ArshonOddSynch} $|x| < 3n$. Therefore, $|x|/q \leq (3n-1)/(n^2+n-1) <
n/(2n-2)$ for all $n \geq 3$, and so $p/q = (|x|+q)/q < (3n-2)/(2n-2) \leq E(\mb a_n)$.
\end{proof}

%\textbf{Note:} For $n = 3$, we get $(3n-1)/(n^2+n-1) = 8/11 > 3/4 = n/(2n-2)$. To prove
%Lemma~\ref{lem:noMordent} for $n=3$, we need to show that actually there are more than $n-1$
%additional $\varphi$-blocks between the blocks containing the mordents. Though this is true, we
%would rather not go into more detailed analysis, since the result for $n=3$ is known.

By Lemma~\ref{lem:noMordent}, in order to compute $E(\mb a_n)$ it is enough to consider powers
$xyx$ such that $x$ is unsynchronized and contains no mordents. The longest subword that contains
no mordents is of length $n+2$, but such subword implies a far pair, and has a unique ancestor.
Therefore, we can assume $|x|\leq n+1$.

\begin{lem}\label{lem:x<=n}
Let $n \geq 3$ be odd, and let $z = xyx = (xy)^{p/q}\in\OC(\mb a_n)$ be an unstretchable power,
such that $x$ is unsynchronized, $x$ contains no mordents, and $|x| = n+1$. Then $p/q < E(\mb
a_n)$.
\end{lem}
\begin{proof}
Since $|x| = n+1$ and $x$ contains no mordents, necessarily $x = ivi$, where $i\in\Sigma_n$ and
either $v = (i+1)\cdots(n-1)01\cdots(i-2)(i-1)$, or $v = (i-1)\cdots01(n-1)\cdots(i+2)(i+1)$.
Suppose the letters of $v$ are increasing, and assume without loss of generality that $i = 0$. Then
$x$ admits two possible interpretations: $x = 01\cdots(n-1)|0$ or $x = 0|1\cdots(n-1)0$. The
ancestors of the first and second case are $\inv{x} = 0^{(e)}1^{(o)}$ and $\inv{x} =
0^{(o)}1^{(e)}$, respectively. Any other interpretation is impossible, since it implies $\mb a_n$
contains two consecutive even $\varphi$-blocks.

As in the previous lemma, we can assume that the two $x$ occurrences of $z$ have different inverse
images. There are two possible cases:
\begin{eqnarray}
\nonumber
\overbrace{01\cdots(n-1)}^{\varphi_{e,n}(0)}|\overbrace{0\cdots}^{\varphi_{o,n}(1)}|&\overbrace{\cdots\cdots\cdots}^{n-2
\;\;\varphi-\textrm{blocks}}&|\overbrace{\cdots
0}^{\varphi_{o,n}(0)}|\overbrace{1\cdots(n-1)0}^{\varphi_{e,n}(1)}\;\;,\\
\nonumber&&\\
\nonumber\overbrace{\cdots 0}^{\varphi_{o,n}(0)}|\overbrace{1\cdots(n-1)0}^{\varphi_{e,n}(1)}|
&\overbrace{\cdots\cdots\cdots}^{n-2 \;\;\varphi-\textrm{blocks}}&|
\overbrace{01\cdots(n-1)}^{\varphi_{e,n}(0)}|\overbrace{0\cdots}^{\varphi_{o,n}(1)}\;\;.
\end{eqnarray}
By Lemma~\ref{lem:dist}, in both cases $y$ contains at least $n-2$ additional $\varphi$-blocks.
Therefore, $q \geq n^2 - n +1$, and so $|x|/q \leq (n+1)/(n^2 - n +1) < n/(2n-2)$ for all $n\geq
3$. Again, $p/q < (3n-2)/(2n-2) \leq E(\mb a_n)$. \end{proof}

By Lemma~\ref{lem:x<=n}, to compute $E(\mb a_n)$ for an odd $n \geq 3$ it is enough to consider
powers of the form $z = xyx$ such that $|x|\leq n$ and $x$ contains no mordent. By
Lemma~\ref{lem:ArshonCE}, such powers have exponent at most $(3n-2)/(2n-2)$. This completes the
proof of Theorem~\ref{thm:ArshonCE}.

\section{Acknowledgement}
I would like to thank Yuri Pritykin, for providing many interesting details on Arshon's life.


\begin{thebibliography}{99}

\bibitem{Arshon:1935} S. E. Arshon. A proof of the existence of infinite asymmetric sequences on $n$
symbols. \emph{Matematicheskoe Prosveshchenie (Mathematical Education)} \textbf{2} (1935), 24--33
(in Russian). Available electronically at \url{http://ilib.mccme.ru/djvu/mp1/mp1-2.htm}.

\bibitem{Arshon:1937} S. E. Arshon. A proof of the existence of infinite asymmetric sequences on $n$ symbols.
\emph{Mat. Sb.} \textbf{2} (1937), 769-–779 (in Russian, with French abstract).

\bibitem{Berstel:1979} J. Berstel. Mots sans carr\'e et morphismes it\'er\'es. \emph{Discrete
Math.} \textbf{29} (1979), 235--244.

\bibitem{Berstel:1995} J. Berstel. Axel Thue's papers on repetitions in words: a translation.
Publications du Laboratoire de Combinatoire et d'Informatique Math\'{e}matique \textbf{20},
Universit\'{e} du Qu\'{e}bec \`{a} Montr\'{e}al (1995).

\bibitem{Currie:2002} J. D. Currie. No iterated morphism generates any Arshon sequence of odd
order. \emph{Discrete Math.} \textbf{259} (2002), 277--283.

\bibitem{Kitaev:2000} S. Kitaev. Symbolic sequences, crucial words and iterations of a morphism.
PhD thesis, G\"oteborg, Sweden, October 2000.

\bibitem{Kitaev:2003} S. Kitaev. There are no iterative morphisms that define the Arshon sequence and the
$\sigma$-sequence. \emph{J. Automata, Languages, and Combinatorics} \textbf{8} no. 1 (2003),
43--50.

\bibitem{Klepinin&Sukhanov:2001} A. V. Klepinin and E. V. Sukhanov. On combinatorial properties of the {Arshon}
sequence. \emph{Disc. Appl. Math.} \textbf{114} (2001), 155--169.

\bibitem{Seebold:2002} P. {S\'{e}\'{e}bold}. About some overlap-free morphisms on a $n$-letter alphabet.
\emph{J. Automata, Languages, and Combinatorics} \textbf{7} (2002), 579--597.

\bibitem{Seebold:2003} P. {S\'{e}\'{e}bold}. On some generalizations of the {Thue--Morse} morphism.
\emph{Theoret. Comput. Sci.} \textbf{292} (2003), 283--298.

\bibitem{Thue:1912} A. Thue, \"{U}ber die gegenseitige Lage gleicher Teile gewisser Zeichenreihen.
\emph{Norske vid. Selsk. Skr. Mat. Nat. Kl.} \textbf{1} (1912), 1--67.

\bibitem {Vilenkin:1991} N. Ya. Vilenkin. Formulas on cardboard.  \emph{Priroda} \textbf{6} (1991),
95--104 (in Russian). English summary available at \url{http://www.ams.org/mathscinet/index.html},
review no. MR1143732.

\end{thebibliography}
\end{document}